\documentstyle[12pt]{amsart}
\begin{document}
\title{Biharmonic maps into compact Lie groups and integrable systems}
\author{Hajime URAKAWA}
\address{Division of Mathematics\\
Graduate School of Information Sciences\\
Tohoku University\\
Aoba 6-3-09, Sendai, 980-8579, Japan
\newline
\quad Current Address: Institute for International Education, Tohoku University\\
Kawauchi 41, Aoba, Sendai, 980-8576, Japan}
\title[Biharmonic maps into compact Lie groups]
{Biharmonic maps into compact Lie groups and integrable systems}
\email{urakawa@@math.is.tohoku.ac.jp}
\thanks{Supported by the Grant-in-Aid for the Scientific Research, 
(A), No. 19204004; (C), No. 21540207, 
Japan Society for the Promotion of Science.}
\dedicatory{}
\maketitle
\begin{abstract} 
In this paper, the reduction 
of biharmonic map equation 
in terms of the Maurer-Cartan form 
for all smooth map 
of a compact Riemannian manifold into a 
compact Lie group $(G,h)$ with the bi-invariant Riemannian metric $h$ is obtained. Due to this formula, 
all the biharmonic curves into compact Lie groups are determined, and 
all the biharmonic maps 
of an open domain of ${\Bbb R}^2$ 
equipped with a Riemannian metric  conformal to 
the standard Euclidean metric into $(G,h)$ are characterized. 
\end{abstract}
\numberwithin{equation}{section}
\theoremstyle{plain}
\newtheorem{df}{Definition}[section]
\newtheorem{th}[df]{Theorem}
\newtheorem{prop}[df]{Proposition}
\newtheorem{lem}[df]{Lemma}
\newtheorem{cor}[df]{Corollary}
\newtheorem{rem}[df]{Remark}
\section{Introduction and statement of results.} 
The theory of harmonic maps of a Riemann surface 
into Lie groups, symmetric spaces 
or homogeneous spaces has been extensively studied 
in the relation with the integrable systems
(\cite{Uh}, \cite{BP}, \cite{W}, \cite{GO}, \cite{DPW}, \cite{DSU}). 
 Let us recall the theory of harmonic maps $\psi$ of a 
 Riemann surface $M$ into a compact Lie group $G$, 
 briefly. 
 Harmonic map is a critical map of the energy functional defined by 
 $$
 E(\psi):=\frac12\int_M\vert d\psi\vert^2\,v_g. 
 $$
 For such a map $\psi$, 
 let $\alpha$ be the pull back of the Maurer-Cartan form $\theta$ of $G$
 which is decomposed into the sum of holomorphic part and the anti-holomorphic one as 
 $\alpha=\alpha'+\alpha''$. 
 Then, it satisfies that $d\alpha+\frac12[\alpha\wedge \alpha]=0$ 
 (the integrability condition),
 and  the harmonicity of $\psi$ is equivalent to 
 the condition $\delta\alpha=0$. 
 Introducing a parameter 
 $\lambda\in {\mathbb C}^{\ast}={\mathbb C}\backslash\{0\}$ 
 as 
 $$
 \alpha_{\lambda}:=
 \frac12(1-\lambda)\alpha'+\frac12(1-\lambda^{-1})\alpha'',
 $$
both the harmonicity and the integrability condition are equivalent to 
 $$
 d\alpha_{\lambda}+\frac12[\alpha_{\lambda}\wedge\alpha_{\lambda}]=0,
 $$
 which implies that there exists an extended solution 
 $\Phi_{\lambda}:\,M\rightarrow G$ satisfying $\Phi_{\lambda}^{-1}d\Phi_{\lambda}=\alpha_{\lambda}$ (\cite{Uh}).  
 Guest and Ohnita (\cite{GO}) showed that 
 the loop group $\Lambda G^{\mathbb C}$ of $G$ acts on the space of all harmonic maps of $M$ into $G$, and Uhlenbeck (\cite{Uh}) showed that 
 every harmonic map from the two-sphere into $G$ is a harmonic map of finite uniton number, and Wood (\cite{W}) determined explicitly 
 harmonic maps of finite uniton numbers. 
 \par
 On the other hand, the theory of biharmonic maps was initiated 
 by Eells and Lemaire (\cite{EL}) and Jiang (\cite{J}). 
 A {\em biharmonic map}
 is a natural extension of harmonic map, and is a critical map of the bienergy functional 
 defined by 
 $$
 E_2(\psi):
 =\frac12\int_M\vert \delta d\psi\vert^2\,v_g=\frac12\int_M\vert \tau(\psi)\vert^2\,v_g,
 $$
 where $\tau(\psi)$ is the tension field of $\psi$, 
 and, by definition, 
 $\psi$ is harmonic if and only if $\tau(\psi)\equiv 0$. 
 \par
 In this paper, we study 
 biharmonic maps of a compact Riemannian manifold $(M,g)$ 
 into a compact Lie group $(G,h)$ with the bi-invariant Riemannian metric $h$. 
 For every $C^{\infty}$ map $\psi: \,(M,g)\rightarrow (G,h)$, 
  let us consider again 
  the pullback $\alpha$ of the Maurer-Cartan form $\theta$. 
  We first will show that the biharmonicity condition for $\psi$ 
  is that 
  $$
  \delta d\delta\alpha+{\rm Trace}_g([\alpha,d\delta \alpha])=0
  $$ 
  (cf. Corollary 3.5) which is a natural extension of harmonicity. 
  Due to this formula, 
  we can determine all real analytic biharmonic curves 
  into a compact Lie group $(G,h)$
 in terms of the initial data 
 $F(0)$, $F'(0)$ and $F''(0)$, 
 where 
 $F(t)=\alpha\left(\frac{\partial}{\partial t}
 \right)$ (cf. Section 4).   
  We give a characterization of biharmonic maps of $({\mathbb R}^2,\mu^2g_0)\rightarrow (G,h)$, 
  where $g_0$ is the standard Euclidean metric on ${\mathbb R}^2$ and 
  $\mu$ is a positive real analytic function on ${\mathbb R}^2$ 
  (cf. Sections 5, 6 and 7). 
\vskip0.6cm\par
{\bf Acknowledgement}: The author expresses his gratitude to Prof. J. Inoguchi 
who gave many useful suggestions and 
Prof. A. Kasue for his financial support 
during the preparation of this paper, and 
Dr. Y. Takenaka and the referees who read carefully and pointed out 
several mistakes in the first draft.  
\vskip0.6cm\par
\section{Preliminaries.}
In this section, we prepare general 
materials and  facts on 
harmonic maps, biharmonic maps into Riemannian manifolds (cf. \cite{EL}, \cite{J}, \cite{KN}). 
\par
Let $(M,g)$ be an $m$-dimensional compact Riemannian manifold, and 
$(N,h)$, an $n$-dimensional Riemannian manifold. 
\par
The {\it energy functional} 
on the space $C^{\infty}(M,N)$ of all $C^{\infty}$ maps 
of $M$ into $N$ is defined by 
$$E(\psi)=\frac12\int_M\vert d\psi\vert^2\,v_g,$$
and 
for a compactly supported 
$C^{\infty}$ one parameter deformation 
$\psi_t\in C^{\infty}(M,N)$ $(-\epsilon<t<\epsilon)$ 
of  $\psi$ with $\psi_0=\psi$, 
the {\it first variation formula} is given by 
$$
\frac{d}{dt}
\bigg\vert _{t=0}E(\psi_t)
=-\int_M
\langle 
\tau(\psi),V \rangle \,v_g,
$$ 
where $V$ is a variation vector field along $\psi$ defined by 
$V=\frac{d}{dt}\big\vert_{t=0} 
\psi_t$ 
which belongs to the space 
$\Gamma(\psi^{-1}TN)$ of sections of the induced bundle of the tangent bundle $TN$ by $\psi$. 
The {\it tension field} $\tau(\psi)$ is 
defined by 
\begin{equation}
\tau(\psi)=-\delta(d\psi),
\end{equation}
where recall 
the definition $\delta \alpha$ 
for a 
$\psi^{-1}TN$-valued $1$-form 
$\alpha$, 
$$
\delta\alpha=-\sum_{i=1}^m
(\overline{\nabla}_{e_i}\alpha)
(e_i)=-\sum_{i=1}^m
\left\{
\overline{\nabla}(\alpha(e_i))
-\alpha(\nabla_{e_i}e_i)
\right\}. 
$$ 
Here, $\nabla$, $\nabla^h$, 
and $\overline{\nabla}$ 
are the Levi-Civita connections of $(M,g)$, $(N,h)$, and the induced connections 
on the induced bundle $\psi^{-1}TN$ 
from $\nabla^h$, respectively. 
For a harmonic map $\psi:\,(M,g)\rightarrow (N,h)$, 
the {\it second variation formula} of the energy functional 
$E(\psi)$ is 
$$
\frac{d^2}{dt^2}\bigg\vert_{t=0}
E(\psi_t)=\int_M\langle J(V),V\rangle\,v_g
$$
where 
\begin{align*}
J(V)&=\overline{\Delta}V-{\mathcal R}(V),\\
\overline{\Delta}V&=\overline{\nabla}^{\ast}\,\overline{\nabla}V
=
-\sum_{i=1}^m
\{
\overline{\nabla}_{e_i}(\overline{\nabla}_{e_i}V)
-\overline{\nabla}_{\nabla_{e_i}e_i}V
\},\\
{\mathcal R}(V)&=\sum_{i=1}^m
R^h(V,d\psi(e_i))d\psi(e_i).
\end{align*}
Here, $\overline{\nabla}$ is the induced connection 
on the induced bundle $\psi^{-1}TN$, and 
is 
$R^h$ is the curvature tensor of $(N,h)$ 
given by 
$R^h(U,V)W=
[\nabla^h_U,\nabla^h_V]W-\nabla^h_{[U,V]}W$ 
$(U,V,W\in {\frak X}(N)$).  
\par
The {\it bienergy functional} is 
defined by 
\begin{equation}
E_2(\psi)=
\frac12
\int_M\vert \delta d\psi\vert^2\,v_g
=\frac12\int_M\vert\tau(\psi)\vert^2\,v_g, 
\end{equation}
and the {\it first variation formula} of the bienergy 
is given (\cite{J}) by 
\begin{equation}
\frac{d}{dt}\bigg\vert _{t=0}E_2(\psi_t)
=-\int_M\langle \tau_2(\psi),V\rangle\,v_g
\end{equation}
where the {\it bitension field} $\tau_2(\psi)$ is 
defined by 
\begin{equation}
\tau_2(\psi)=J(\tau(\psi))
=\overline{\Delta}\tau(\psi)-{\mathcal R}(\tau(\psi)),
\end{equation}
and a $C^{\infty}$ map
$\psi:(M,g)\rightarrow (N,h)$ is called to be 
{\it biharmonic} if 
\begin{equation}
\tau_2(\psi)=0.
\end{equation}
\par
The biharmonic maps are real analytic when both $(M,g)$ and $(N,h)$ are real analytic. Because the solutions of non-linear elliptic partial differential equations are real analytic. 
\vskip0.6cm\par
\section{Determination of the bitension field}
Now, 
assume that  $(N,h)$ is 
an $n$-dimensional compact Lie group 
with Lie algebra $\frak g$, 
and $h$, 
the bi-invariant Riemannian metric on $G$ 
corresponding to the Ad$(G)$-invariant inner product 
$\langle\,,\,\rangle$ on $\frak g$. 
Let $\theta$ be the Maurer-Cartan form on $G$, i.e., 
a $\frak g$-valued left invariant $1$-form on $G$ 
which is 
defined by 
$\theta_y(Z_y)=Z$, 
($y\in G$, 
$Z\in \frak g$). 
For every 
$C^{\infty}$ map $\psi$ of $(M,g)$ into $(G,h)$, 
let us consider 
a $\frak g$-valued $1$-form $\alpha$ on $M$ 
given by 
$\alpha=\psi^{\ast}\theta$. 
Then it is well known (see for example, 
\cite{DSU}) that 
\begin{lem} 
For every $C^{\infty}$ map 
$\psi:\,(M,g)\rightarrow (G,h)$, 
\begin{equation}
\theta(\tau(\psi))=-\delta\alpha.
\end{equation}
\par
Thus, $\psi:(M,g)\rightarrow (G,h)$ is harmonic if and only if $\delta\alpha=0$. 
\end{lem}
\vskip0.6cm\par
Let 
$\{X_s\}_{s=1}^n$ be an orthonormal basis 
of $\frak g$ with respect to 
the inner product $\langle\,,\,\rangle$. 
Then, for every 
$V\in \Gamma(\psi^{-1}TG)$, 
\begin{align}
&V(x)=\sum_{s=1}^nh_{\psi(x)}(V(x),X_{s\,\psi(x)})\,X_{s\,\psi(x)}\in T_{\psi(x)}G,\nonumber\\
&\theta(V)(x)=\sum_{s=1}^nh_{\psi(x)}(V(x),X_{s\,\psi(x)})\,X_s\in \frak g,
\end{align}
for all $x\in M$. 
Then, for every $X\in {\frak X}(M)$, 
\begin{align}
\theta(\overline{\nabla}_XV)&=
\sum_{s=1}^nh(\overline{\nabla}_XV,X_{s})\,X_s\nonumber\\
&
=\sum_{s=1}^n\{
X\,h(V,X_s)-h(V,\overline{\nabla}_XX_s)
\}X_s\nonumber\\
&=X(\theta(V))-\sum_{s=1}^nh(V,\overline{\nabla}_XX_s)X_s,
\end{align}
where we regarded a vector field $Y\in {\frak X}(G)$ by 
$Y(x)=Y(\psi(x)) \,(x\in M)$
to be an element in the space 
$\Gamma(\psi^{-1}TG)$ of smooth sections of $\psi^{-1}TG$.
\par
Here, let us recall that the Levi-Civita connection 
$\nabla^h$ of $(G,h)$ is given 
(cf. \cite{KN} Vol. II, p. 201, Theorem 3.3)
by 
\begin{equation}
\nabla^h_{X_t}X_s=\frac12\,[X_t,X_s]
=\frac12\sum_{\ell=1}^nC^{\ell}_{ts}\,X_{\ell}, 
\end{equation}
where the structure constant 
$C_{ts}^{\ell}$ 
of $\frak g$ 
is defined by 
$
[X_t,X_s]=\sum_{\ell=1}^nC_{ts}^{\ell}\,X_{\ell}, 
$
and satisfies that 
\begin{align}
C^{\ell}_{ts}
=\langle[X_t,X_s],X_{\ell}\rangle
=-\langle X_s,[X_t,X_{\ell}]\rangle
=-C_{t\ell}^s.
\end{align}
Thus, we have by (3.4) and (3.5), 
\begin{align}
\sum_{s=1}^nh(V,\overline{\nabla}_XX_s)X_s
&=\frac12
\sum_{s,t=1}^n
h\left(V,\sum_{\ell=1}^nh(\psi_{\ast}X,X_t)C^{\ell}_{ts}\,X_{\ell}\right)\,X_s\nonumber\\
&=-\frac12\sum_{s,t,\ell=1}^nh(V,X_{\ell})\,
h(\psi_{\ast}X,X_t)\,C^s_{t\ell}X_s\nonumber\\
&=-\frac12\sum_{t,\ell=1}^nh(V,X_{\ell})\,h(\psi_{\ast}X,X_t)\,[X_t,X_{\ell}]\nonumber\\
&=-\frac12\,\left[
\sum_{t=1}^nh(\psi_{\ast}X,X_t)\,X_t,
\sum_{\ell=1}^nh(V,X_{\ell})\,X_{\ell}
\right]\nonumber\\
&=-\frac12\,\left[
\alpha(X),\theta(V)
\right],
\end{align}
because
we have 
\begin{align}
\alpha(X)
=\theta(\psi_{\ast}X)
=\sum_{t=1}^nh(\psi_{\ast}X,X_t)X_t,
\end{align}
and 
\begin{align}
\theta(V)=\sum_{\ell=1}^nh(V,X_{\ell})\,\theta(X_{\ell})
=\sum_{\ell=1}^nh(V,X_{\ell})\,X_{\ell}.
\end{align}
\par
Therefore, inserting (3.6) into (3.3), 
we obtain 
\begin{lem} 
For every $C^{\infty}$ map $\psi:\,(M,g)\rightarrow (G,h)$, 
\begin{equation}
\theta(\overline{\nabla}_XV)
=X(\theta(V))+\frac12\,[\alpha(X),\theta(V)], 
\end{equation}
where 
$V\in \Gamma(\psi^{-1}TG)$ and 
$X\in {\frak X}(M)$. 
\end{lem}
\par\qed
\vskip0.6cm\par
We shall show 
\begin{th}
For every $\psi\in C^{\infty}(M,G)$, we have 
\begin{align}
\theta(\tau_2(\psi))&=\theta(J(\tau(\psi)))\nonumber\\
&=-\delta\,d\,\delta\alpha-{\rm Trace}_g([\alpha,d\,\delta\alpha]),
\end{align}
where $\alpha=\psi^{\ast}\theta$. 
\end{th}
\vskip0.6cm\par
Here, let us recall the definition:
\begin{df}
For two $\frak g$-valued $1$-forms 
$\alpha$ and $\beta$ on $M$, we define 
$\frak g$-valued symmetric $2$-tensor 
$[\alpha,\beta]$ on $M$ 
by 
\begin{equation}
[\alpha,\beta](X,Y):=\frac12
\left\{
[\alpha(X),\beta(Y)]+[\alpha(Y),\beta(X)]
\right\},\quad (X,Y\in{\frak X}(M))
\end{equation}
and its trace ${\rm Trace}_g([\alpha,\beta])$ by 
\begin{equation}
{\rm Trace}_g([\alpha,\beta])
:=\sum_{i=1}^m[\alpha,\beta](e_i,e_i).
\end{equation}
\par 
Recall a $\frak g$-valued $2$-form 
$[\alpha\wedge \beta]$ on $M$ 
in the introduction to be 
given by 
\begin{equation}
[\alpha\wedge \beta](X,Y)
:=\frac12\left\{
[\alpha(X),\beta(Y)]
-[\alpha(Y),\beta(X)]\right\},
\quad(X,Y\in {\frak X}(M)).
\end{equation}
\end{df}
Then, we have immediately by Theorem 3.3, 
\begin{cor}
For every $\psi\in C^{\infty}(M,G)$, we have 
\par
$(1)$ $\psi:\,(M,g)\rightarrow (G,h)$ is harmonic if and only if 
\begin{equation}\delta\alpha=0.\end{equation} 
\par
$(2)$ $\psi:\,(M,g)\rightarrow (G,h)$ is biharmonic if and only if 
\begin{equation}\delta\,d\,\delta\alpha+{\rm Trace}_g([\alpha,d\,\delta\alpha])=0.\end{equation} 
\end{cor}
\vskip0.8cm\par
We give a proof of Theorem 3.3.
\begin{pf}
(The first step) \quad 
We first show that, for all $V\in \Gamma(\psi^{-1}TG)$, 
\begin{align}
\theta(\overline{\Delta}V)&=
\Delta_g\theta(V)\nonumber
-\sum_{i=1}^m
\left\{
\frac12[e_i(\alpha(e_i)),\theta(V)]
+[\alpha(e_i),e_i(\theta(V))]
\right.\nonumber\\
&
\qquad
+\frac14[\alpha(e_i),[\alpha(e_i),\theta(V)]]
\left.
-\frac12[\alpha(\nabla_{e_i}e_i),\theta(V)]
\right\},
\end{align}
where 
$\{e_i\}_{i=1}^m$ is a locally defined 
orthonormal frame field on $(M,g)$, and 
$\Delta_g$ is the (positive) Laplacian 
of $(M,g)$ acting on $C^{\infty}(M)$. 
\par
Indeed, we have by using Lemma 3.2 twice, 
\begin{align}
\theta(\overline{\Delta}V)
&=-\sum_{i=1}^m
\left\{
\theta(\overline{\nabla}_{e_i}(\overline{\nabla}_{e_i}V))
-\theta(\overline{\nabla}_{\nabla_{e_i}e_i}V)
\right\}\nonumber\\
&=
-\sum_{i=1}^m
\left\{
e_i(\theta(\overline{\nabla}_{e_i}V))
+\frac12[\alpha(e_i),\theta(\overline{\nabla}_{e_i}V)]
\right.
\nonumber\\
&\left.\qquad\qquad
-\nabla_{e_i}e_i(\theta(V))
-\frac12[\alpha(\nabla_{e_i}e_i),\theta(V)]\right\}\nonumber\\
&=
-\sum_{i=1}^m
\left\{
e_i\left(e_i(\theta(V))+\frac12[\alpha(e_i),\theta(V)]
\right)
\right.\nonumber\\
&\qquad\qquad
+\frac12
\left[
\alpha(e_i),e_i(\theta(V))+\frac12[\alpha(e_i),\theta(V)]
\right]\nonumber\\
&\qquad\qquad
\left.
-\nabla_{e_i}e_i\,(\theta(V))-\frac12[\alpha(\nabla_{e_i}e_i),\theta(V)]
\right\}\nonumber\\
&=
-\sum_{i=1}^m\left\{
e_i(e_i(\theta(V)))-\nabla_{e_i}e_i\,(\theta(V))\right\}\nonumber\\
&\quad
-\sum_{i=1}^m\left\{
\frac12
e_i([\alpha(e_i),\theta(V)])
+\frac12[\alpha(e_i),e_i(\theta(V))]
\right.\nonumber\\
&
\qquad\quad
+\frac14[\alpha(e_i),[\alpha(e_i),\theta(V)]]
\left.
-\frac12[\alpha(\nabla_{e_i}e_i),\theta(V)]
\right\}.
\end{align}
Here,  we have 
\begin{equation*}
e_i([\alpha(e_i),\theta(V)])
=[e_i(\alpha(e_i)),\theta(V)]+[\alpha(e_i),e_i(\theta(V))],
\end{equation*}
which we substitute into (3.17), and
by definition of $\Delta_g$, 
we have (3.16). 
\vskip0.3cm\par
(The second step) \quad 
On the other hand, 
we have to consider 
\begin{align}
-\sum_{i=1}^mR^h(V,\psi_{\ast}e_i)\psi_{\ast}e_i
=-\sum_{i=1}^mR^h(L_{\psi(x)\,\ast}^{-1}V,
L_{\psi(x)\,\ast}^{-1}\psi_{\ast}e_i)L_{\psi(x)\,\ast}^{-1}\psi_{\ast}e_i.
\end{align}
Under the identification 
$T_eG\ni Z_e\leftrightarrow Z\in {\frak g}$, 
we have 
\begin{align}
&T_eG\ni L_{\psi(x)\,\ast}^{-1}\psi_{\ast}e_i\leftrightarrow \alpha(e_i)\in {\frak g},\\
&T_eG\ni L_{\psi(x)\,\ast}^{-1}V\leftrightarrow \theta(V)\in {\frak g},
\end{align}
respectively. 
Because, we have 
$$
L_{\psi(x)\,\ast}^{-1}\psi_{\ast}e_i
=\sum_{s=1}^n
h(\psi_{\ast}e_i,X_{s\,\psi(x)})\,X_{s\,e}
$$
and 
\begin{align*}
\alpha(e_i)&=\psi^{\ast}\theta(e_i)
=\theta(\psi_{\ast}e_i)
=\sum_{s=1}^nh(\psi_{\ast}e_i,X_{s\,\psi(x)})\,
\theta(X_{s\,\psi(x)})\\
&=\sum_{s=1}^nh(\psi_{\ast}e_i,X_{s\,\psi(x)})\,X_s,
\end{align*}
which implies that (3.19). Analogously, we obtain (3.20). 
\vskip0.3cm\par
Under this identification, the curvature tensor of 
$(G,h)$ is given as (see 
Kobayashi-Nomizu (\cite{KN}, pp. 203--204)),  
$$R^h(X,Y)_e=-\frac14{\rm ad}([X,Y])
\qquad(X,Y\in {\frak g}),
$$
and then, we have 
\begin{align}
\theta\left(-\sum_{i=1}^mR^h(V,\psi_{\ast}e_i))\psi_{\ast}e_i\right)
&=\frac14\sum_{i=1}^m
\left[
\left[
\theta(V),\alpha(e_i)
\right],\alpha(e_i)
\right]\nonumber\\
&=\frac14\sum_{i=1}^m
\left[
\alpha(e_i),\left[
\alpha(e_i),\theta(V)\right]
\right].
\end{align}
\vskip0.3cm\par
(The third step) \quad 
By (3.16) and (3.21), 
for $V\in \Gamma(\psi^{-1}TG)$, we have 
\begin{align}
\theta
&\left(\overline{\Delta}V
-\sum_{i=1}^m
R^h(V,\psi_{\ast}e_i)\psi_{\ast}e_i\right)\nonumber\\
&=
\Delta_g\theta(V)\nonumber\\
&
\quad
-\sum_{i=1}^m
\left\{
\frac12[e_i(\alpha(e_i)),\theta(V)]
+[\alpha(e_i),e_i(\theta(V))]
\right.
+\frac14[\alpha(e_i),[\alpha(e_i),\theta(V)]]\nonumber\\
&\qquad\qquad\qquad\quad
\left.
-\frac12[\alpha(\nabla_{e_i}e_i),\theta(V)]
\right\}\nonumber\\
&
\quad
+\frac14
\sum_{i=1}^m[\alpha(e_i),[\alpha(e_i),\theta(V)]]\nonumber\\
&=
\Delta_g\theta(V)-\frac12\sum_{i=1}^me_i(\alpha(e_i)),\theta(V)]+\sum_{i=1}^m[\alpha(e_i),e_i(\theta(V))]
\nonumber\\
&\qquad\qquad\quad
+\frac12\sum_{i=1}^m[\alpha(\nabla_{e_i}e_i),\theta(V)]\nonumber\\
&=
\Delta_g\theta(V)
-\frac12\left[
\sum_{i=1}^m\left(
e_i(\alpha(e_i))-\alpha(\nabla_{e_i}e_i)
\right),\theta(V)
\right]
+\sum_{i=1}^m[\alpha(e_i),e_i(\theta(V))]\nonumber\\
&=
\Delta_g\theta(V)
+\frac12[\delta\alpha,\theta(V)]
+\sum_{i=1}^m[\alpha(e_i),e_i(\theta(V))].
\end{align}
\vskip0.3cm\par
(The fourth step) \quad For $V=\tau(\psi)$ in (3.22), 
since  
$\theta(\tau(\psi))=-\delta\alpha$, we have 
\begin{align}
\theta(J(\tau(\psi)))&=
\Delta_g\theta(\tau(\psi))
+\frac12[\delta\alpha,\theta(\tau(\psi))]
\nonumber\\
&\qquad\qquad+\sum_{i=1}^m[\alpha(e_i),e_i(\theta(\tau(\psi))]\nonumber\\
&=-\Delta_g\delta\alpha-\frac12[\delta\alpha,\delta\alpha]
-\sum_{i=1}^m[\alpha(e_i),e_i(\delta\alpha)]\nonumber\\
&=-\Delta_g\delta\alpha-\sum_{i=1}^m[\alpha(e_i),e_i(\delta\alpha)]\nonumber\\
&=\Delta_g\delta\alpha
-\sum_{i=1}^m[\alpha(e_i),d\delta\alpha(e_i)].
\end{align}
Then, (3.23) implies the desired (3.10).
\end{pf}
\vskip0.6cm\par
\section{Biharmonic curves from ${\mathbb R}$ into compact Lie groups}
In this section, we consider the simplest case: 
$(M,g)=({\mathbb R},g_0)$ is 
the standard $1$-dimensional Euclidean space, 
 and $(G,h)$ is an $n$-dimensional compact Lie group with the bi-invariant Riemannian metric $h$. 
\vskip0.6cm\par
{\bf 4.1} \quad 
First, let 
$\psi:\,{\mathbb R}\ni t\mapsto \psi(t)\in 
(G,h)$, a $C^{\infty}$ curve in $G$. Then, $\alpha:=\psi^{\ast}\theta$ is a $\frak g$-valued $1$-form on 
$\mathbb R$. So, $\alpha$ can be written 
at $t\in {\mathbb R}$ as 
\begin{equation}
\alpha_t=F(t)\,dt
\end{equation}
where $F:\,{\mathbb R}\ni t\mapsto F(t)\in  {\frak g}$ 
is a $C^{\infty}$ function on ${\mathbb R}$ given by 
\begin{equation}
F(t)=\alpha\left(
\frac{\partial}{\partial t}
\right)
=\psi^{\ast}\theta
\left(
\frac{\partial}{\partial t}
\right)
=\theta\left(
\psi_{\ast}\left(
\frac{\partial}{\partial t}
\right)
\right).
\end{equation}
Here, 
since 
\begin{equation}
\psi'(t):=\psi_{\ast}\left(
\frac{\partial}{\partial t}
\right)
=\sum_{s=1}^n
h_{\psi(t)}\left(
\psi_{\ast}\left(
\frac{\partial}{\partial t}
\right), 
X_{s\,\psi(t)}
\right)\,X_{s\,\psi(t)},
\end{equation}
we have
\begin{equation}
F(t)=\sum_{s=1}^n
h_{\psi(t)}
\left(
\psi_{\ast}
\left(
\frac{\partial}{\partial t}
\right),X_{s\,\psi(t)}
\right)\,X_s,
\end{equation}
so that we have the following correspondence: 
\begin{align}
T_eG\ni
&L_{\psi(t)\,\ast}^{-1}\psi'(t)
=\sum_{s=1}^nh_{\psi(t)}(\psi'(t),X_{s\,\psi(t)})\,X_{s\,e}
\nonumber\\
&\leftrightarrow 
F(t)=
\theta
\left(
\psi_{\ast}\left(
\frac{\partial}{\partial t}
\right)
\right)
\in {\frak g}. 
\end{align}
\vskip0.6cm\par
{\bf 4.2} \quad 
We have  that 
\begin{equation}
\delta\alpha=-F'(t), 
\end{equation}
since we have 
$
\delta\alpha
=-e_1(\alpha(e_1))
=-e_1(F(t))=-F'(t).
$
\vskip0.3cm\par
Therefore, we have
$\psi:\,({\mathbb R},g_0)\rightarrow (G,h)$ 
is {\it harmonic} if and only if
\begin{align}
\delta\alpha=0
\quad&\Longleftrightarrow
\quad 
F'=0\nonumber\\
\quad&\Longleftrightarrow\quad 
\alpha=X\otimes dt
\quad(\text{for some} \,\,X\in {\frak g})\nonumber\\
\quad&\Longleftrightarrow\quad 
\psi:\,{\mathbb R}\rightarrow \,(G,h),\,\text{\it geodesic},
\end{align}
since 
\begin{equation}
F(t)=\theta(\psi'(t))=L_{\psi(t)\ast}{}^{-1}\psi'(t),
\end{equation}
we have 
\begin{equation}
\psi'(t)=L_{\psi(t)\ast}X=X_{\psi(t)},
\end{equation}
for some $X\in {\frak g}$ 
which yields that 
$$
\psi(t)=x\,\exp(tX).
$$
\par
Therefore, {\it any geodesic through} $\psi(0)=x$ 
{\it is given by }
\begin{equation}
\psi(t)=x\,\exp(tX),\,\quad(t\in {\mathbb R})\,\,
\end{equation}
\text{\it for some} 
$X\in {\frak g}$. 
\par
On the other hand, we want to determine 
a {\it biharmonic curve} 
$\psi:\,({\mathbb R},g_0)\rightarrow (G,h)$.  
By (4.6), we have 
\begin{equation}
\delta d\delta\alpha
=-\frac{\partial^2}{\partial t^2}
\left(
-F'(t)
\right)
=F^{(3)}(t), 
\end{equation}
and 
\begin{align}
{\rm Trace}_g[\alpha,d\delta\alpha]
=\left[
\alpha\left(
\frac{\partial}{\partial t}
\right),\,d \delta\alpha\left(
\frac{\partial}{\partial t}
\right)
\right]
=\left[
F(t),F''(t)
\right],
\end{align}
so by (4.9), (4.10), and (3.16) in Corollary 3.5, 
$\psi:\,({\mathbb R},g_0)\rightarrow (G,h)$ is 
{\it biharmonic} if and only if 
\begin{equation}
F^{(3)}-[F(t),F''(t)]=0.
\end{equation}
\par
 \vskip0.6cm\par
{\bf 4.3} \quad 
For a $C^{\infty}$ curve $\psi:\,{\mathbb R}\rightarrow G$, 
let $\psi(t):=\exp X(t)$, where $X(t)\in {\frak g}$. 
Then, 
\begin{equation}
F(t)=\theta
\left(
\psi_{\ast}
\left(
\frac{\partial}{\partial t}
\right)\right),
\,\,\psi_{\ast}
\left(
\frac{\partial}{\partial t}
\right)
\in T_{\psi(t)}G,
\end{equation}
and  by the following formula (cf. \cite{H}, p. 95) 
\begin{equation*}
\exp_{\ast\,X}
=L_{\exp\,X\,\ast\,e}
\circ\,
\frac{1-e^{-\,{\rm ad\,X}}}{{\rm ad}\,X}
\quad (X\in {\frak g}), 
\end{equation*}
we have 
\begin{align} 
\psi_{\ast}
\left(
\frac{\partial}{\partial t}
\right)
&=\exp_{\ast\,X(t)}X'(t) 
\nonumber\\
&=L_{\exp \,X(t)\,\ast \,e}
\left(
\sum_{n=0}^{\infty}
\frac{(- {\rm ad}\,X(t))^n}{(n+1)!}
\,(X'(t))
\right).
\end{align}
Since $\theta$ is a left invariant $1$-form, 
 we have 
 \begin{equation}
 F(t)=\sum_{n=0}^{\infty}
\frac{(- {\rm ad}\,X(t))^n}{(n+1)!}
\,(X'(t)). 
 \end{equation}
 {\bf 4.4} \quad 
 The initial value problem 
 \begin{equation}
 \left\{
 \begin{aligned}
 &F^{(3)}(t)=[F(t),F''(t)],\\
 &F(0)=B_0,F'(0)=B_1,F''(0)=B_2,
 \end{aligned}
 \right.
 \end{equation}
 for every $B_i\in {\frak g}$ $(i=0,1,2)$, 
 has a unique solution $F(t)$. 
 Assume that 
 $X(t)$ is a real analytic curve in $t$, and $X(0)=0$. 
 Then, 
 $F(t)$ is also real analytic in $t$, and 
 we can write as 
 \begin{equation}
 X(t)=\sum_{n=1}^{\infty}A_n\,t^n,\quad 
 F(t)=\sum_{n=0}^{\infty}B_n\,t^n. 
 \end{equation}
 By (4.16), 
 we have 
 \begin{align}
 F(t)&=X'(t)+\frac12[-X(t),X'(t)]+\frac16[-X(t),[-X(t),X'(t)]]\nonumber\\
 &\quad +\sum_{n=3}^{\infty}
 \frac{(-{\rm ad}\,X(t))^n}{(n+1)!}\,(X'(t)). 
 \end{align}
Since 
 $X'(t)=\sum_{m=0}^{\infty}A_{m+1}(m+1)\,t^m$, 
 we have 
 $$\frac12[-X(t),X'(t)]=-\frac12[A_1,A_2]\,t^2+O(t^3),$$
 and 
 $$\frac16[-X(t),[-X(t),X'(t)]]=O(t^3),$$
 so that we have 
 $$
 F(t) = A_1+2A_2\,t+\left(3A_3-\frac12[A_1,A_2]\right)t^2+O(t^3).
 $$
 We continue this process, we have 
 \begin{equation}
 \left\{
 \begin{aligned}
 B_0&=A_1\\
 B_1&=2A_2,\\
 B_2&=3A_3-\frac12[A_1,A_2]\\
 \cdots&\cdots\cdots\cdots\cdots\cdots\cdots\cdots\cdots\\
 B_n&=(n+1)A_{n+1}+G_n(A_1,\cdots,A_n),
 \end{aligned}
 \right.
 \end{equation}
 where $G_n(x_1,\cdots,x_n)$ is a polynomial in $(x_1,\cdots,x_n)$. 
 Notice that for arbitrary given data $(B_0,B_1,B_2)$, 
 all $B_n$ $(n=0,1,\cdots)$ are determined, and by using (4.20), 
 one can determine all $A_n$ $(n=1,2,\cdots)$, uniquely. 
 Therefore, by summarizing the above, we obtain 
 \begin{th}
 For every $C^{\infty}$ curve $\psi:\,{\mathbb R}\rightarrow G$, 
 $\psi(t)=\exp \,X(t)$, $(X(t)\in {\frak g})$ and 
 \begin{equation}
 \alpha\left(\frac{\partial}{\partial t}\right)=F(t)=\sum_{n=0}^{\infty}
\frac{(- {\rm ad}\,X(t))^n}{(n+1)!}
\,(X'(t)). 
 \end{equation}
 \par
 $(1)$ $\psi:\,({\mathbb R},g_0)\rightarrow (G,h)$ is biharmonic 
 if and only if 
 \begin{equation}
 F^{(3)}(t)=[F(t),F''(t)].
 \end{equation}
 \par
 $(2)$ 
 The initial value problem 
 \begin{equation}
 \left\{
 \begin{aligned}
 &F^{(3)}(t)=[F(t),F''(t)],\\
 &F(0)=B_0,F'(0)=B_1,F''(0)=B_2,
 \end{aligned}
 \right.
 \end{equation}
 has a unique solution $F(t)$ for arbitrary given data 
 $(B_0,B_1,B_2)$ in ${\frak g}$. 
 \par
 $(3)$ Assume that $\psi:\,({\mathbb R},g_0)\rightarrow (G,h)$ 
 is a real analytic biharmonic curve with $\psi(0)=e$. 
 Then, 
 $\psi(t)$ is uniquely determined by 
 $F(0)=B_0$, $F'(0)=B_1$, and $F''(0)=B_2$.  
 \end{th}
 \vskip0.6cm\par
 {\it Example} 
 \quad If $G$ is abelian, let us consider 
 a $C^{\infty}$ curve $\psi:\,{\mathbb R}\rightarrow G$ given by  
 $\psi(t)=\exp\,X(t)$. 
 Then, $F(t)=X'(t)$, 
 and 
 $\psi:\,({\mathbb R},g_0)\rightarrow (G,h)$ 
 is biharmonic if and only if 
 $F^{(3)}(t)=X^{(4)}(t)=0$. 
 Then, $X(t)=A_0+A_1\,t+A_2\,t^2+A_3\,t^3$. 
 Thus, every biharmonic curve 
 $\psi:\,({\mathbb R},g_0)\rightarrow (G,h)$ with $\psi(0)=e$ is given by
 $$
 \psi(t)=\exp(A_1\,t+A_2\,t^2+A_3\,t^3). 
 $$ 
 \vskip0.6cm\par
 {\bf 4.5} \quad
 Now we will solve the ODE $(4.22)$ for 
 a biharmonic isometric immersion 
 $\psi:\,({\mathbb R},g_0)\rightarrow G$ and  
 a $\frak g$-valued 
 curve $F(t)$ in the case of ${\frak g}={\frak su}(2)$. 
  Let  $G=SU(2)$ with the 
 bi-invariant Riemannian metric $h$ 
 which corresponds to the following 
 Ad$(SU(2))$-invariant inner product 
 $\langle\,,\,\rangle$ 
 on 
 \begin{align*}
 &{\frak g}=
 {\frak s}{\frak u}(2)
 =\{X\in M(2,{\mathbb C});\,X+{}^{\rm t}\overline{X}=0,
 {\rm Tr}(X)=0\}.\\
& \langle X,Y\rangle=
 -2{\rm Tr}(XY)\quad (X,Y\in {\frak s}{\frak u}(2)).
 \end{align*}
 If we choose 
 $$
 X_1=\begin{pmatrix}
 \frac{\sqrt{-1}}{2}&0\\
 0&-\frac{\sqrt{-1}}{2}
 \end{pmatrix},\,\,
  X_2=\begin{pmatrix}
 0&\frac12\\
 -\frac12&0
 \end{pmatrix},\,\,
 X_3=\begin{pmatrix}
 0&\frac{\sqrt{-1}}{2}\\
 \frac{\sqrt{-1}}{2}&0
 \end{pmatrix}, 
 $$
 then $\{X_1,X_2,X_3\}$ is an orthonormal basis of 
 $({\frak s}{\frak u}(2),\langle\,,\,\rangle)$, and satisfies 
 the Lie bracket relations: 
 $$
 [X_1,X_2]=X_3,\,\,[X_2,X_3]=X_1,\,\,[X_3,X_1]=X_2.
 $$
 Thus, the ODE's (4.22) becomes 
 \begin{equation}
 \left\{
 \begin{aligned}
 y_1^{(3)}&=y_2\,y_3''-y_3\,y_2'',\\
 y_2^{(3)}&=y_3\,y_1''-y_1\,y_3'',\\
 y_3^{(3)}&=y_1\,y_2''-y_2\,y_1'',
 \end{aligned}
 \right.
 \end{equation}
 which is equivalent to 
 \begin{equation}
 {\bf y}^{(3)}={\bf y}\times {\bf y}'',
 \end{equation}
 where 
 ${\bf y}:={}^{\rm t}\!(y_1,y_2,y_3)\in {\mathbb R}^3$, and ${\bf a}\times{\bf b}$ 
 stands for the vector cross product in 
 ${\mathbb R}^3$. 
 Notice here that even so $\frak g$ is non-abelian, 
 but our equation (4.22) is the vector equation depending 
 on the time $t$ of the Euclidean space 
 ${\mathbb R}^3$ by identifying 
 $\sum_{i=1}^3y_i\,X_i\in {\frak g}\mapsto (y_1,y_2,y_3)\in {\mathbb R}^3$. 
 \par
 Then, the ODE's $(4.25)$ can be solved as follows: 
 \par
 Let ${\bf x}(s)={}^{\rm t}\!(x_1(s),x_2(s),x_3(s))$ be a $C^{\infty}$ curve in ${\mathbb R}^3$ with arc length parameter $s$, and then 
 \begin{align*}
 {\bf y}(s)&={\bf x}'(s)={\bf e}_1(s).
  \end{align*}
  Let $\{{\bf e}_1(s),{\bf e}_2(s),{\bf e}_3(s)\}$ 
  be the Frenet frame field along 
  ${\bf x}(s)$. Recall the Frenet-Serret formula:
 \begin{equation*}
 \left\{
 \begin{aligned}
 {\bf e}_1{}'&=\qquad \quad\kappa \,{\bf e}_2\\
 {\bf e}_2{}'&=-\kappa\,{\bf e}_1\qquad\quad+\tau\,{\bf e}_3\\
 {\bf e}_3{}'&=\qquad\quad -\tau\,{\bf e}_2
 \end{aligned}
 \right.
 \end{equation*}
 where $\kappa$ and $\tau$ are the curvature and torsion of 
 ${\bf x}(s)$, respectively. 
 Then, we have 
 \begin{equation}
 \left\{
 \begin{aligned}
 {\bf y}'&=\kappa\,{\bf e}_2\\
 {\bf y}''&=-\kappa^2\,{\bf e}_1+\kappa'\,{\bf e}_2+\kappa\tau\,{\bf e}_3\\
 {\bf y}'''&=
 -3\kappa\kappa'\,{\bf e}_1+(\kappa''-\kappa^3-\kappa\tau^2)\,{\bf e}_2+(2\kappa'\tau+\kappa\tau')\,{\bf e}_3.
 \end{aligned}
 \right.
 \end{equation}
 Thus, (4.24) is equivalent to 
 \begin{align}
 -3\kappa\kappa'\,{\bf e}_1
 &+(\kappa''-\kappa^3-\kappa\tau^2)\,{\bf e}_2+(2\kappa'\tau+\kappa\tau')\,{\bf e}_3\nonumber\\
 &=
 {\bf e}_1\times
 (-\kappa^2\,{\bf e}_1+\kappa'\,{\bf e}_2
 +\kappa\tau\,{\bf e}_3)
 \nonumber\\
 &=-\kappa\tau\,{\bf e}_2+\kappa'\,{\bf e}_3
 \end{align}
 which is equivalent to 
 \begin{equation}
 \left\{
 \begin{aligned}
 -3\kappa\kappa'&=0\\
 \kappa''-\kappa^3-\kappa\tau^2&=-\kappa\tau\\
 2\kappa'\tau+\kappa\tau'&=\kappa'.
 \end{aligned}
 \right.
 \end{equation}
 \par
 Then, the first equation of (4.28) 
 turns out that 
 $(\kappa^2)'=0$, that is, 
 $\kappa ^2$ is constant, i.e., 
 $\kappa\equiv 0$, or 
 $\kappa\equiv \kappa_0\not=0$. 
 In the case that $\kappa\equiv0$, 
 the solution of (4.28),  
 ${\bf x}(s)$, is a line in ${\mathbb R}^3$. 
 \par
 For the case that 
 $\kappa\equiv\kappa_0\not=0$, 
 the only solution of (4.24) is 
 \begin{equation}
 \left\{
 \begin{aligned}
 \kappa&\equiv\kappa_0\not=0,\\
 \tau&\equiv\tau_0, \,\,\text{and}\\
 \kappa_0{}^2&=\tau_0(1-\tau_0),
 \end{aligned}
 \right.
 \end{equation}
and the unique solution of (4.25) is  given by 
 \begin{equation}
 {\bf x}(s)=
 \begin{pmatrix}
 x_1(s)\\
 x_2(s)\\
 x_3(s)
 \end{pmatrix}
 =\begin{pmatrix}
 a\,\cos\frac{s}{\sqrt{a^2+1}}+b\\
 a\,\sin\frac{s}{\sqrt{a^2+1}}+b\\
 \frac{s}{\sqrt{a^2+1}}+b
 \end{pmatrix}
 \end{equation}
 for some positive constant $a>0$ and 
 some constant $b$. 
 Thus, $F(s)$  
 is given as follows:
 \begin{align}
 F(s)&={\bf x}'(s)=\sum_{i=1}^3{x_i}'(s)\,X_i \nonumber\\
 &=
 \left(
 -\frac{a}{\sqrt{a^2+1}}\,\sin\frac{s}{\sqrt{a^2+1}}
 \right)X_1
 +\left(
 \frac{a}{\sqrt{a^2+1}}\,\cos\frac{s}{\sqrt{a^2+1}}
 \right)X_2
 \nonumber\\
 &\quad
 +
 \left(
 \frac{1}{\sqrt{a^2+1}}
 \right)X_3,
 \end{align}
 for any constant $a>0$. 
 Conversely, it is easy to see that every such $F(s)$ in $(4.31)$ 
 is a solution of 
 $(4.22)$: $F^{(3)}(s)=[F(s),F''(s)]$. 
 \vskip0.6cm\par
 {\it Remark} \quad 
 It is still difficult for us to determine 
 $X(t)$ to satisfy $(4.21)$: 
 $$
 F(t) =\sum_{n=0}^{\infty}
 \frac{(-{\rm ad} X(t))^n}{(n+1)!}(X'(t)),
 $$
 in the case of ${\frak su}(2)$.  
 \vskip0.6cm\par
 \section{Biharmonic maps from an open domain 
 in ${\mathbb R}^2$}
 In this section, we consider a biharmonic map 
 $\psi:\,({\mathbb R}^2,g)\supset \Omega
 \rightarrow (G,h)$. 
 Here, we assume that 
 $G$ is a linear compact Lie group, i.e., 
 $G$ is a subgroup of the unitary group 
 $U(N)(\subset GL(N,{\mathbb C}))$ of degree $N$  
 with a bi-invariant Riemannian metric $h$ on $G$. 
 Let $\frak g$ be the Lie algebra of $G$ which is 
  a Lie subalgebra of 
  the Lie algebra
  ${\frak u}(N)$ of
  $U(N)$. 
 The Riemannian metric $g$ on ${\mathbb R}^2$ 
 is a conformal metric which is given by 
 $g=\mu^2\,g_0$ with 
 a $C^{\infty}$ positive function $\mu$ on $\Omega$ and 
 $g_0=dx\cdot dx+dy\cdot dy$, 
 where  
 $(x,y)$ is the standard coordinate on ${\mathbb R}^2$. 
 \par
 Let 
 $\psi:\,\Omega\ni (x,y)\mapsto 
 \psi(x,y)=(\psi_{ij}(x,y))\in U(N)$ a $C^{\infty}$ map. 
  Let us consider 
  $$
  \frac{\partial \psi}{\partial x}:=\left(
  \frac{\partial \psi_{ij}}{\partial x}
  \right), 
  \quad 
  \frac{\partial \psi}{\partial y}:=\left(
  \frac{\partial \psi_{ij}}{\partial y}
  \right). 
  $$
  Then, 
  \begin{equation}
  A_x:= \psi^{-1}\frac{\partial\psi}{\partial x},
  \quad 
  A_y:= \psi^{-1}\frac{\partial\psi}{\partial y}
  \end{equation}
  are ${\frak g}$-valued 
  $C^{\infty}$ functions on $\Omega$. 
  It is known that, 
  for two given $\frak g$-valued 
  $1$-forms $A_x$ and $A_y$ on $\Omega$,
  there exists a $C^{\infty}$ mapping 
  $\psi:\,\Omega\rightarrow G$ 
  satisfying the equations 
  (5.1) 
  if the {\it integrability condition} holds:
  \begin{equation}
  \frac{\partial A_y}{\partial x}
  - \frac{\partial A_x}{\partial y}
  +[A_x,A_y]=0.
  \end{equation}
  The pull back of the Maurer-Cartan form $\theta$ 
  by $\psi$ is given by 
  \begin{align}
  \alpha&:=\psi^{\ast}\theta=
  \psi^{-1}d\psi
  =\psi^{-1}\frac{\partial\psi}{\partial x}\,dx
  +\psi^{-1}\frac{\partial\psi}{\partial y}\,dy\nonumber\\
  &=A_x\,dx+A_y\,dy,
  \end{align}
  which is a $\frak g$-valued 
  $1$-form on $\Omega$. 
  \par
  Recall that 
  the codifferential $\delta\alpha$ is 
  of a $\frak g$-valued $1$-form 
  $\alpha=A_x\,dx+A_y\,dy$, 
  $A_x=\psi^{-1}\frac{\partial\psi}{\partial x}$ 
  and 
  $A_y=\psi^{-1}\frac{\partial\psi}{\partial y}$, is given by 
  \begin{equation}
  \delta\alpha=-
  \mu^{-2}
  \left\{
  \frac{\partial}{\partial x}A_x+\frac{\partial}{\partial y}A_y
  \right\}.
  \end{equation}
  Then, we have the following well known facts: 
  \begin{lem}
  We have 
  \begin{align}
  \delta\alpha
  &=-\mu^{-2}
  \left\{
  \frac{\partial}{\partial x}
  \left(
  \psi^{-1}\frac{\partial\psi}{\partial x}
  \right)
  +
  \frac{\partial}{\partial x}
  \left(
  \psi^{-1}\frac{\partial\psi}{\partial y}
  \right)
  \right\}\\
  &=-\mu^{-2}
  \left\{
  \frac{\partial A_x}{\partial x}+
  \frac{\partial A_y}{\partial y}
  \right\}. 
  \end{align} 
 Therefore, the following three statements are equivalent: 
 \begin{align}
  &(i) \qquad\qquad 
  \psi:\,(\Omega,g)\rightarrow (G,h)
  \,\,
  \text{is {\em harmonic}},  \nonumber
  \\
 &(ii) \qquad\qquad
  \delta\alpha=0,\\
  &(iii) \qquad \qquad
 \frac{\partial A_x}{\partial x}+
  \frac{\partial A_y}{\partial y}=0.
  \end{align}
  \end{lem}
  \vskip0.6cm\par
  Next, calculate 
  the Laplacian $\Delta_g$ of $({\mathbb R}^2,g)$
  for $g=\mu^2\,g_0$. 
  We obtain 
  \begin{align}
  \Delta_g&=-\sum_{i,j=1}^2
  g^{ij}\left(
  \frac{\partial^2}{\partial x^i\,\partial x^j}-\sum_{k=1}^2\Gamma^k_{ij}\,\frac{\partial}{\partial x^k}
  \right)\nonumber\\
  &=
  -\mu^{-2}\left(
  \frac{\partial^2}{\partial x^2}+
  \frac{\partial^2}{\partial y^2}
  \right).
  \end{align}
  Thus we have 
  \begin{align}
  \delta d\delta\alpha
  &=
  \Delta_g(\delta\alpha)\nonumber\\
  &=
  \mu^{-2}
  \,\left(
  \frac{\partial^2}{\partial x^2}+
   \frac{\partial^2}{\partial y^2}
  \right)
  \left[
  \mu^{-2}
  \left\{
  \frac{\partial}{\partial x}
  \left(
  \psi^{-1}\frac{\partial \psi}{\partial x}
  \right)
  +\frac{\partial}{\partial y}
  \left(
  \psi^{-1}\frac{\partial \psi}{\partial y}
  \right)
  \right\}
  \right]\nonumber\\
  &=\mu^{-2}
  \,\left(
  \frac{\partial^2}{\partial x^2}+
   \frac{\partial^2}{\partial y^2}
  \right)
  \left[
  \mu^{-2}
  \left\{
  \frac{\partial A_x}{\partial x}
  +\frac{\partial A_y}{\partial y}
  \right\}
  \right]\nonumber\\
  &=-\mu^{-2}
  \left(
  \frac{\partial^2}{\partial x^2}+
   \frac{\partial^2}{\partial y^2}
  \right)
(\delta\alpha).
  \end{align}
  \par
  On the other hand, 
  by taking 
  an orthonormal local frame field 
  $\{e_1,e_2\}$ of $({\mathbb R}^2,g)$, as
  $
  e_1=\mu^{-1}\frac{\partial}{\partial x},
  \quad 
  e_2=\mu^{-1}\frac{\partial}{\partial y},
  $
  we have 
  \begin{align}
  {\rm Trace}_g&([\alpha,d\delta\alpha])
  =[\alpha(e_1),d\delta\alpha(e_1)]
  +[\alpha(e_2),d\delta\alpha(e_2)]\nonumber\\
  &=
 - \mu^{-2}\left[
  A_x,
  \frac{\partial}{\partial x}
 \left( \mu^{-2}
  \left\{
  \frac{\partial A_x}{\partial x}
  +\frac{\partial A_y}{\partial y}
  \right\}\right)
  \right]\nonumber\\
  &\quad
  -\mu^{-2}\left[
  A_y,
  \frac{\partial}{\partial y}
  \left(\mu^{-2}
  \left\{
  \frac{\partial A_x}{\partial x}
  +\frac{\partial A_y}{\partial y}
  \right\}\right)
  \right]\nonumber\\
  &=
  \mu^{-2}[A_x,\frac{\partial}{\partial x}(\delta\alpha)]
  +
   \mu^{-2}[A_y,\frac{\partial}{\partial y}(\delta\alpha)].
  \end{align}
  By (5.10) and (5.11), we obtain 
  \begin{align}
  \delta d\delta\alpha&+{\rm Trace}_g([\alpha,d\delta\alpha])\nonumber\\
  &=
  -\mu^{-2}
  \left(
  \frac{\partial^2}{\partial x^2}+
   \frac{\partial^2}{\partial y^2}
  \right)(\delta\alpha)
 +
 \mu^{-2}[A_x,\frac{\partial}{\partial x}(\delta\alpha)]
  +
   \mu^{-2}[A_y,\frac{\partial}{\partial y}(\delta\alpha)] 
  \nonumber\\
  &=
  -\mu^{-2}
  \left\{
   \left(
  \frac{\partial^2}{\partial x^2}+
   \frac{\partial^2}{\partial y^2}
  \right)(\delta\alpha)\right.
  \left.
  -\frac{\partial}{\partial x}[A_x,\delta\alpha]
  -\frac{\partial}{\partial y}[A_y,\delta\alpha]
  \right\},
  \end{align}
  where in the last equation in (5.11), 
  we only notice that 
  \begin{align*}
   \frac{\partial}{\partial x}[A_x,&\delta\alpha]
  +\frac{\partial}{\partial y}[A_y,\delta\alpha]\nonumber\\
&=
\left[\frac{\partial}{\partial x}A_x,\delta\alpha\right]
+\left[A_x,\frac{\partial}{\partial x}(\delta\alpha)\right]
+\left[\frac{\partial}{\partial y}A_y,\delta\alpha\right]
+\left[A_y,\frac{\partial}{\partial y}(\delta\alpha)\right]\nonumber\\
&=
\left[
\frac{\partial}{\partial x}A_x+
\frac{\partial}{\partial y}A_y, \delta\alpha
\right]
+\left[A_x,\frac{\partial}{\partial x}(\delta\alpha)\right]
+\left[A_y,\frac{\partial}{\partial y}(\delta\alpha)\right]
\nonumber\\
&=[-\mu^{-2}\delta\alpha,\delta\alpha]
+\left[A_x,\frac{\partial}{\partial x}(\delta\alpha)\right]
+\left[A_y,\frac{\partial}{\partial y}(\delta\alpha)\right]\nonumber\\
&=\left[A_x,\frac{\partial}{\partial x}(\delta\alpha)\right]
+\left[A_y,\frac{\partial}{\partial y}(\delta\alpha)\right].
  \end{align*}
  Thus, we have
  \begin{th}
  Let $\Omega$ be an open subset of ${\mathbb R}^2$,  $g=\mu^{2}g_0$,
  a Riemannian metric
  conformal to
 the standard metric $g_0$ on 
 $\Omega$ with a $C^{\infty}$ positive function $\mu$ on $\Omega$, and 
 $\psi:\,\Omega\rightarrow G$, 
 a $C^{\infty}$ map of $\Omega$ into  
 a compact linear Lie group $(G,h)$
 with bi-invariant Riemannian metric $h$. 
  Then, 
  \par
 $(1)$  The $1$-form $\alpha$ satisfies 
   $d\alpha+\frac12[\alpha\wedge\alpha]=0$ which is equivalent to
   \begin{equation}
   \frac{\partial A_y}{\partial x}-\frac{\partial A_x}{\partial y}+[A_x,A_y]=0. 
   \end{equation}
  \par
    $(2)$ The following three are equivalent: 
  \begin{align}
  &(i)
  \qquad
  \psi:\,(\Omega,g)\rightarrow (G,h)
 \text{  is {\em harmonic}}, \nonumber\\
 &(ii)
 \qquad  \delta\alpha=0,\\
& (iii)\qquad
  \frac{\partial}{\partial x}A_x
  +\frac{\partial}{\partial y}A_y=0. 
   \end{align}
   \par
  $(3)$ The following three are equivalent: 
  \begin{align} 
  &(i) \qquad
  \psi:\,(\Omega,g)\rightarrow (G,h)
  \text{  is {\em biharmonic}},\nonumber\\
  &(ii) \quad\,\,\,
  \delta d \delta\alpha+{\rm Trace}_g([\alpha,d\delta\alpha])=0,\\
  &(iii) \quad 
   \left(
  \frac{\partial^2}{\partial x^2}+
   \frac{\partial^2}{\partial y^2}
  \right)(\delta\alpha)
  -\frac{\partial}{\partial x}[A_x,\delta\alpha]
  -\frac{\partial}{\partial y}[A_y,\delta\alpha]
=0.
  \end{align}
  \par
  $(4)$ 
  Let us consider two $\frak g$-valued $1$-forms 
  $\beta$ and $\Theta$ on 
  $\Omega$, defined by
  \begin{align}
  \beta&:=[A_x,\delta\alpha]\,dx+[A_y,\delta\alpha]\,dy, \\
  \Theta&:=d\delta\alpha-\beta,
  \end{align}
  respectively. 
  Then, $\psi:(\Omega,g)\rightarrow (G,h)$ is 
  {\em biharmonic} if and only if 
  \begin{equation} 
  \delta\Theta=0.
   \end{equation}
   \end{th}
  \begin{pf} 
  (1) is clear. We see already (2) and (3). 
  For $(4)$, we only have to see that 
 (5.17) is equivalent to 
 \begin{align}
 0=-\Delta_g(\delta\alpha)+\delta\beta=-\delta(d\delta\alpha-\beta)=-\delta\Theta
 \end{align} 
 where
 \begin{align}
 \Theta&:=
 d\delta\alpha-\beta\nonumber\\
 &=\frac{\partial}{\partial x}(\delta\alpha)\,dx
 +\frac{\partial}{\partial y}(\delta\alpha)\,dy
 -[A_x,\delta\alpha]\,dx
 -[A_y,\delta\alpha]\,dy\nonumber\\
 &=
 \left\{
 \frac{\partial}{\partial x}(\delta\alpha)
 -[A_x,\delta\alpha]
 \right\}\,dx
 +
  \left\{
 \frac{\partial}{\partial y}(\delta\alpha)
 -[A_y,\delta\alpha]
 \right\}\,dy.
 \end{align}
  \end{pf}
  \vskip0.6cm\par
  \section{Complexfication of the biharmonic map equation}
  We use the complex coordinate 
  $z=x+iy$ $(i=\sqrt{-1})$ in $\Omega$, and we put 
  $A_z=\frac12(A_x-iA_y)$ and 
  $A_{\overline{z}}=\frac12(A_x+iA_y)$ 
  which are ${\frak g}^{\mathbb C}$-valued functions with 
  $A_{\overline{z}}=\overline{A_z}$.  
  Then, it is well known that
  \begin{align*}
  \frac{\partial}{\partial\overline{z}}A_z
  +\frac{\partial}{\partial z}A_{\overline{z}}
  &=\frac12\left\{
  \frac{\partial}{\partial x}A_x+
   \frac{\partial}{\partial y}A_y
  \right\},\\
  \frac{\partial}{\partial z}A_{\overline{z}}
  -\frac{\partial}{\partial\overline{z}}A_z
  &+[A_z,A_{\overline{z}}]
  =
  \frac{i}{2}
  \left\{
  \frac{\partial}{\partial x}A_y
  - \frac{\partial}{\partial y}A_x+[A_x,A_y]
  \right\}, 
  \end{align*}
 and also
 \begin{align*}
&\alpha=A_x\,dx+A_y\,dy
 =A_z\,dz+A_{\overline{z}}\,d\overline{z},\\
 &\frac{\partial^2}{\partial x^2}
 + \frac{\partial^2}{\partial y^2}
 =4\,\frac{\partial^2}{\partial z\partial\overline{z}},\\
 &\delta\alpha=
 -\mu^{-2}
 \left(
 \frac{\partial}{\partial x}A_x
 + \frac{\partial}{\partial y}A_y
 \right)
 =
 -2\mu^{-2}\left(
 \frac{\partial}{\partial\overline{z}}A_z
 +\frac{\partial}{\partial z}A_{\overline{z}}
 \right).
\end{align*}
Then, the conditions (5.20) is 
equivalent to that 
\begin{equation}
\delta\widetilde{\Theta}=0,
\end{equation}
where
\begin{equation}
\widetilde{\Theta}:=
 \left\{
 \frac{\partial}{\partial z}(\delta\alpha)
 -[A_z,\delta\alpha]
 \right\} dz
 +
  \left\{
 \frac{\partial}{\partial \overline{z}}(\delta\alpha)
 -[A_{\overline{z}},\delta\alpha]
 \right\}d\overline{z}.
\end{equation}
The integrability condition $(5.13)$ is equivalent to 
\begin{equation}
 \frac{\partial}{\partial z}A_{\overline{z}}
  -\frac{\partial}{\partial\overline{z}}A_z
  +[A_z,A_{\overline{z}}]=0. 
\end{equation}
\section{Determination of biharmonic maps}
In this section, we want to show how to determine all the biharmonic maps of 
$(\Omega,g)$ into a compact Lie group $(G,h)$ 
where $g=\mu^{2}g_0$ with a positive $C^{\infty}$  function on $\Omega$ and 
$h$ is a bi-invariant Riemannian metric on $G$. 
\par
Our method to obtain all the biharmonic maps is divide three steps: 
\par
({\it The first step}) \quad 
We solve first the equation:
\begin{equation}
\frac{\partial}{\partial\overline{z}}
B_z+\frac{\partial}{\partial z}B_{\overline{z}}=0.
\end{equation}
Notice that, if these $B_z$ and $B_{\overline{z}}$ satisfy furthermore, 
the integrability condition 
\begin{equation}
\frac{\partial}{\partial z}B_{\overline{z}}
-\frac{\partial}{\partial \overline{z}}B_z
+[B_z,B_{\overline{z}}]=0,
\end{equation}
then, there exists a harmonic map 
$\Psi:\,(\Omega,g)\rightarrow(G,h)$ such that 
\begin{equation}
\left\{
\begin{aligned}
\Psi^{-1}\frac{\partial\Psi}{\partial z}&=B_z,\\
\Psi^{-1}\frac{\partial\Psi}{\partial \overline{z}}&
=B_{\overline{z}}, 
\end{aligned}
\right.
\end{equation}
and the converse is true. 
\par
({\it The second step}) \quad 
For such 
two ${\frak g}^{\mathbb C}$-valued functions $B_z$ and 
$B_{\overline{z}}$ on $\Omega$ satisfying $(7.1)$ not necessarily satisfying $(7.2)$, 
we should detect two 
${\frak g}^{\mathbb C}$-valued functions $A_z$ and 
$A_{\overline{z}}$ on $\Omega$ satisfying that 
\begin{equation}
\left\{
\begin{aligned}
&\frac{\partial}{\partial z}
\left(
-2\mu^{-2}
\left(
\frac{\partial A_z}{\partial\overline{z}}
+\frac{\partial A_{\overline{z}}}{\partial z}
\right)
\right)
-\left[A_z,
-2\mu^{-2}
\left(
\frac{\partial A_z}{\partial\overline{z}}
+\frac{\partial A_{\overline{z}}}{\partial z}
\right)
\right]
=B_z,\\
&\frac{\partial}{\partial \overline{z}}
\left(
-2\mu^{-2}
\left(
\frac{\partial A_z}{\partial\overline{z}}
+\frac{\partial A_{\overline{z}}}{\partial z}
\right)
\right)-\left[
A_{\overline{z}},-2\mu^{-2}\left(
\frac{\partial A_z}{\partial\overline{z}}
+\frac{\partial A_{\overline{z}}}{\partial z}
\right)\right]=B_{\overline{z}},\\
&\frac{\partial}{\partial z}
A_{\overline{z}}-\frac{\partial}{\partial\overline{z}}A_z
+[A_z,A_{\overline{z}}]=0.
\end{aligned}
\right.
\end{equation}
\par
({\it The third step})\quad 
Finally, for the above ${\frak g}^{\mathbb C}$-valued functions
$A_z$ and $A_{\overline{z}}$ on $\Omega$ 
satisfying $(7.4)$  
and $a\in G$, 
there exists a $C^{\infty}$ mapping 
$\psi:\,\Omega\rightarrow G$ satisfying that 
\begin{equation}
\left\{
\begin{aligned}
\psi(x_0,y_0)&=a,\\
\psi^{-1}\frac{\partial \psi}{\partial z}&=A_z,\\
\psi^{-1}\frac{\partial \psi}{\partial \overline{z}}&
=A_{\overline{z}}.
\end{aligned}
\right.
\end{equation}
\par
Then, $\psi:\,(\Omega,g)\rightarrow (G,h)$ is 
a {\it biharmonic} map 
due to $(5.20)$, $(6.1)$ and $(7.4)$, and conversely, 
every biharmonic map $\psi:\,(\Omega,g)\rightarrow (G,h)$ could be obtained in this way. 
\par
To do the these procedures rigorously, 
let us define 
\begin{df}
(1) Let us define the four sets 
$\Lambda$, $\Lambda_1$, $\Lambda_2$, and $\Lambda_0$: 
\newline
$\bullet$ Let 
$\Lambda$ be 
the set of all $\frak g$-valued two functions $(A_x,A_y)$ on $\Omega$, 
(or all ${\frak g}^{\mathbb C}$-valued two functions 
$(A_z,A_{\overline{z}})$ on $\Omega$ 
with $A_{\overline{z}}=\overline{A_z}$, 
\newline$\bullet$ 
let $\Lambda_1$, 
the set of $(A_x,A_y)\in \Lambda$ which satisfy the harmonic map equation 
$(5.12)$ (or  $(7.1)$), 
\newline$\bullet$ 
let $\Lambda_2$, 
the set of $(A_x,A_y)\in \Lambda$ which satisfy the biharmonic map equation $(5.17)$ (or $(6.1)$), and 
\newline$\bullet$ 
let 
$\Lambda_0$, 
the set of $(A_x,A_y)\in\Lambda$ which 
satisfy the integrability condition 
$(5.13)$, (or $(6.3)$), 
respectively. 
\vskip0.3cm\par
(2) Let us define two sets $\Xi$ and $\Xi_1$:
\newline$\bullet$ 
Let $\Xi$ be 
the set of all $\frak g$-valued two real analytic functions $(B_x,B_y)$ 
on $\Omega$ 
(or ${\frak g}^{\mathbb C}$-valued two real analytic functions 
$(B_z,B_{\overline{z}})$ 
on $\Omega$ with $B_{\overline{z}}=\overline{B_z}$), 
and 
\newline$\bullet$ 
let $\Xi_1$,  
the set of all $(B_x,B_y)=(B_z,B_{\overline{z}})\in \Xi$  
satisfying the harmonic map equation $(7.1)$, 
respectively.  
\end{df}
\vskip0.6cm\par
\begin{df}
Let us define two $C^{\infty}$ mappings 
$\Phi_i$ $(i=1,2)$ of $\Lambda$ into $\Xi$ 
by 
\begin{align}
\Phi_1(A_x,A_y)&:=
\bigg(\frac{\partial}{\partial x}\left(-\mu^{-2}\bigg(\frac{\partial A_x}{\partial x}+\frac{\partial A_y}{\partial y}\bigg)\right)
-\left[
A_x,
-\mu^{-2}\bigg(\frac{\partial A_x}{\partial x}+\frac{\partial A_y}{\partial y}\bigg)
\right],
\nonumber\\
&
\frac{\partial}{\partial y}\left(-\mu^{-2}\bigg(\frac{\partial A_x}{\partial x}+\frac{\partial A_y}{\partial y}\bigg)\right)
-\left[
A_y,
-\mu^{-2}\bigg(\frac{\partial A_x}{\partial x}+\frac{\partial A_y}{\partial y}\bigg)
\right]\bigg),
\\
\Phi_2(A_x,A_y)&:=
\bigg(
-\mu^{-2}\left(
\frac{\partial^2A_x}{\partial x^2}+\frac{\partial^2A_x}{\partial y^2}
-\frac{\partial}{\partial y}[A_x,A_y]
\right)\nonumber\\
&-\frac{\partial \mu^{-2}}{\partial x}\,\left(
\frac{\partial A_x}{\partial x}+\frac{\partial A_y}{\partial y}
\right)
-\left[
A_x,-\mu^{-2}\bigg(\frac{\partial A_x}{\partial x}+\frac{\partial A_y}{\partial y}\bigg)
\right],\nonumber\\
&
-\mu^{-2}\left(
\frac{\partial^2A_y}{\partial x^2}+\frac{\partial^2A_y}{\partial y^2}
-\frac{\partial}{\partial x}[A_x,A_y]
\right)\nonumber\\
&-\frac{\partial \mu^{-2}}{\partial y}\,\left(
\frac{\partial A_x}{\partial x}+\frac{\partial A_y}{\partial y}
\right)
-\left[
A_y,-\mu^{-2}\bigg(\frac{\partial A_x}{\partial x}+\frac{\partial A_y}{\partial y}\bigg)
\right]\bigg),
\end{align}
respectively. 
\end{df}
Then, we obtain 
\begin{th} Assume that $\Omega$ be a simply connected open domain 
in ${\mathbb R}^2$, and 
$\mu$ is a positive real analytic function on $\Omega$.  Then, we have: 
\par
$(1)$ For every 
$(B_x,B_y)=(B_z,B_{\overline{z}})\in \Xi$ 
there exists $(A_x,A_y)=(A_z,A_{\overline{z}})\in \Lambda$ 
such that 
$\Phi_2(A_x,A_y)=(B_x,B_y)$ (or $\Phi_2(A_z,A_{\overline{z}})=(B_z,B_{\overline{z}})$). The solution $(A_x,A_y)=(A_z,A_{\overline{z}})$ is uniquely determined 
by the initial data $A_x(x_0,y)$, $A_y(x_0,y)$, 
$\frac{\partial A_x}{\partial x}(x_0,y)$ and $\frac{\partial A_y}{\partial x}(x_0,y)$, $(x_0,y)\in \Omega$. 
\par
$(2)$ $\Phi_1=\Phi_2$ on $\Lambda_0$,
\par
$(3)$ $\Phi_1{}^{-1}(\Xi_1)=\Lambda_2$, and 
$\Phi_1(\Lambda_2\cap \Lambda_0)
=\Phi_2{}(\Lambda_2\cap \Lambda_0)=\Xi_1$. 
\end{th}
\begin{pf} 
 For (1), by definition of $\Phi_2$, 
 that $\Phi_2(A_x,A_y)=(B_x,B_y)$ 
 is equivalent to the following two equations: 
 \begin{align}
 \frac{\partial^2A_x}{\partial x^2}
 &=-\frac{\partial^2A_x}{\partial y^2}
 +\frac{\partial}{\partial y}[A_x,A_y]
\nonumber \\
 &
 \qquad
 -\mu^2\frac{\partial \mu^{-2}}{\partial x}\,
 \left(
 \frac{\partial A_x}{\partial x}+\frac{\partial A_y}{\partial y}
 \right)
 -\mu^2\bigg[
 A_x,-\mu^{-2}
  \left(
 \frac{\partial A_x}{\partial x}+\frac{\partial A_y}{\partial y}
 \right)\bigg]\nonumber\\
 &\qquad -\mu^2\,B_x,\\
  \frac{\partial^2A_y}{\partial x^2}
  &=-\frac{\partial^2A_y}{\partial y^2}
 +\frac{\partial}{\partial x}[A_x,A_y]
\nonumber \\
 &
 \qquad
 -\mu^2\frac{\partial \mu^{-2}}{\partial y}\,
 \left(
 \frac{\partial A_x}{\partial x}+\frac{\partial A_y}{\partial y}
 \right)
 -\mu^2\bigg[
 A_y,-\mu^{-2}
  \left(
 \frac{\partial A_x}{\partial x}+\frac{\partial A_y}{\partial y}
 \right)\bigg]\nonumber\\
 &\qquad -\mu^2\,B_y.
 \end{align}
 Notice that the system of (7.8) and (7.9) satisfies 
 all the conditions of the theorem of Cauchy-Kovalevskaya 
 when $n_i=2$ $(i=1,2)$ (cf. 
 \cite{E}, p. 1305, 429 B; \cite{M}, p. 224; \cite{I}, p. 181)
 \begin{th} (Cauchy-Kovalevskaya) 
 Let us consider the following Cauchy problem of unknown $N$ functions 
 $u_i(t,x)$ $(i=1,\cdots, N)$ 
 in $t$ and $x=(x_1,\cdots,x_m)$,  
 \begin{equation}
 \left\{
 \begin{aligned}
 &\frac{\partial^{n_i}u_i}{\partial t^{n_i}}=F_i(t,x,D^k_tD^p_xu_j)\quad(i=1,\cdots,N),\\
& \frac{\partial^ku_i}{\partial t^k}(t_0,x)=\varphi^k_i(x)
 \,\,(0\leq k\leq n_i-1; \,i=1,\cdots,N),
 \end{aligned}
 \right.
 \end{equation}
 where, 
 for $p=(p_1,\cdots,p_m)$, 
 $\vert p\vert=p_1+\cdots +p_m$, 
 $
 D^k_tD^p_x
 :=\frac{\partial^k}{\partial t^k}\,
 \frac{\partial^{\vert p\vert}}{\partial x_1{}^{p_1}\cdots\partial_m{}^{p_m}}
 $
 and in the right hand side of the first equation of (7.10), 
 $k$ and $p$ satisfy 
 $$
 k<n_j \quad\text{and}\quad k+\vert p\vert\leq n_j \quad (j=1,\cdots,N).
 $$
  Assume that each $F_i$ and $\varphi^k_i$ are real analytic functions. 
 Then, there exists a real analytic solution $u_i$ $(i=1,\cdots,N)$ of 
 (7.10)
 and it is unique in the class of real analytic functions. 
  \end{th}
 \vskip0.6cm\par
 Then, for each $(B_x,B_y)\in \Xi$, 
 there exists a real analytic solution $(A_x,A_y)$ of the Cauchy problem 
 (7.8) and (7.9) with the initial condition:  
 \begin{equation}
 \left\{
 \begin{aligned}
 \left(\frac{\partial A_x}{\partial x}\right)(x_0,y)&=f_1(y),\quad A_x(x_0,y)=f_0(y),\\
  \left(\frac{\partial A_y}{\partial x}\right)(x_0,y)&=g_1(y),\quad A_x(x_0,y)=g_0(y),
 \end{aligned}
 \right.
 \end{equation} 
  and the real analytic solution 
  $(A_x,A_y)$ is unique for real analytic functions 
  $f_i$ and $g_i$ $(i=1,2)$. 
 \par 
 By taking this process at each point $(x_0,y_0)$ in $\Omega$, we have a real analytic 
 solution $(A_x,A_y)$ of (7.8) and (7.9) in an open neighborhood of $(x_0,y_0)$. 
 Then, by the uniqueness theorem of the continuation of a real analytic function on a simply connected domain $\Omega$, 
 we have a solution $(A_x,A_y)$ of (7.8) and (7.9) on $\Omega$. 
 We have $(1)$. 
 \par
 For (2), 
we have to see $\Phi_1(A_x,A_y)=\Phi_2(A_x,A_y)$ 
for every $(A_x,A_y)\in \Lambda_0$, 
which follows from 
that 
\begin{align*}
\frac{\partial}{\partial x}\left(
\mu^{-2}\left(
\frac{\partial A_x}{\partial x}+\frac{\partial A_y}{\partial y}
\right)
\right)&=
\mu^{-2}\left(
\frac{\partial^2 A_x}{\partial x^2}
+\frac{\partial^2A_y}{\partial x\partial y}
\right)\\
&\qquad+\frac{\partial\mu^{-2}}{\partial x}\,\left(
\frac{\partial A_x}{\partial x}+\frac{\partial A_y}{\partial y}
\right)\\
&=
\mu^{-2}\left(
\frac{\partial^2 A_x}{\partial x^2}
+\frac{\partial^2A_x}{\partial y^2}
-\frac{\partial}{\partial y}[A_x,A_y]
\right)\\
&\qquad+\frac{\partial\mu^{-2}}{\partial x}\,\left(
\frac{\partial A_x}{\partial x}+\frac{\partial A_y}{\partial y}
\right),
\end{align*}
because of (5.13) and 
it is a similar for 
$\frac{\partial}{\partial y}\left(
\mu^{-2}\left(
\frac{\partial A_x}{\partial x}+\frac{\partial A_y}{\partial y}
\right)
\right)$, so that we have (2). 
\par
For (3), due to (2), 
we only have to see $\Phi_1{}^{-1}(\Xi_1)=\Lambda_2$ 
which is equivalent to that:
\begin{align*}
\qquad 
&\text{\em for all} \,\,(B_x,B_y)\in \Xi_1,\,\,\text{\em exists a unique}
\,\,(A_x,A_y)\in \Lambda_2
\,\,\text{\em such that}\,\,
\\
\qquad 
&\Phi_1(A_x,A_y)=(B_x,B_y), 
\text{ 26\em and vice versa}. 
\end{align*}
\par 
But, that 
$(B_x,B_y)=(B_z,B_{\overline{z}})\in \Xi_1$ means that 
it satisfies the harmonic map equation (7.1). 
On the other hand, 
$\Phi_1(A_x,A_y)=(B_x,B_y)$ means that 
$\Phi_1(A_z,A_{\overline{z}})
=(B_z,B_{\overline{z}})$ 
which is equivalent to that 
the first two equations of (7.4) hold 
by definition of $\Phi_1$,  
and notice here that 
$\Phi_1(A_x,A_y)=(B_x,B_y)$ is equivalent to 
the two following equations 
\begin{align}
 &\frac{\partial}{\partial x}\left(-\mu^{-2}\bigg(\frac{\partial A_x}{\partial x}+\frac{\partial A_y}{\partial y}\bigg)\right)
-\left[
A_x
-\mu^{-2}\bigg(\frac{\partial A_x}{\partial x}+\frac{\partial A_y}{\partial y}\bigg)
\right] = B_x,
\\
&
\frac{\partial}{\partial y}\left(-\mu^{-2}\bigg(\frac{\partial A_x}{\partial x}+\frac{\partial A_y}{\partial y}\bigg)\right)
-\left[
A_y,
-\mu^{-2}\bigg(\frac{\partial A_x}{\partial x}+\frac{\partial A_y}{\partial y}\bigg)
\right]
=B_y,
\end{align}
which are also equivalent to 
\begin{align}
&\frac{\partial}{\partial z}
\left(
-2\mu^{-2}
\left(
\frac{\partial A_z}{\partial\overline{z}}
+\frac{\partial A_{\overline{z}}}{\partial z}
\right)
\right)
-\left[A_z,
-2\mu^{-2}
\left(
\frac{\partial A_z}{\partial\overline{z}}
+\frac{\partial A_{\overline{z}}}{\partial z}
\right)
\right]
=B_z,\\
&\frac{\partial}{\partial \overline{z}}
\left(
-2\mu^{-2}
\left(
\frac{\partial A_z}{\partial\overline{z}}
+\frac{\partial A_{\overline{z}}}{\partial z}
\right)
\right)-\left[
A_{\overline{z}},-2\mu^{-2}\left(
\frac{\partial A_z}{\partial\overline{z}}
+\frac{\partial A_{\overline{z}}}{\partial z}
\right)\right]=B_{\overline{z}}. 
\end{align}
But, by inserting both 
(7.14) and (7.15) into 
\begin{equation*}
\frac{\partial}{\partial\overline{z}}
B_z+\frac{\partial}{\partial z}B_{\overline{z}}=0.
\qquad \qquad\qquad (7.1), 
\end{equation*}
we obtain 
\begin{align}
&
\frac{\partial^2}{\partial\overline{z}\partial z}
\left(
-2\mu^{-2}
\left(
\frac{\partial A_z}{\partial\overline{z}}
+\frac{\partial A_{\overline{z}}}{\partial z}
\right)
\right)
-\frac{\partial}{\partial\overline{z}}
\left[A_z,
-2\mu^{-2}
\left(
\frac{\partial A_z}{\partial\overline{z}}
+\frac{\partial A_{\overline{z}}}{\partial z}
\right)
\right]\nonumber\\
&
+\frac{\partial^2}{\partial z\partial \overline{z}}
\left(
-2\mu^{-2}
\left(
\frac{\partial A_z}{\partial\overline{z}}
+\frac{\partial A_{\overline{z}}}{\partial z}
\right)
\right)-
\frac{\partial}{\partial z}\left[
A_{\overline{z}},-2\mu^{-2}\left(
\frac{\partial A_z}{\partial\overline{z}}
+\frac{\partial A_{\overline{z}}}{\partial z}
\right)\right]\nonumber\\
&\qquad =0, 
\end{align}
which is just the biharmonic map equation for $(A_z,A_{\overline{z}})$: (6.1) 
$\delta\widetilde{\Theta}=0$. 
By the same way, one can see also immediately 
$(A_x,A_y)$ satisfies 
the biharmonic map equation (5.20) if 
$(B_x,B_y)$ satisfies the harmonic map 
equation (5.15) 
by using Theorem 5.2, 
(5.6) and 
(5.22). 
Thus, we obtain 
$\Phi_1{}^{-1}(\Xi_1)=\Lambda_2$ 
and (3). 
\end{pf}
\vskip0.6cm\par
{\it Remark.} \quad 
 The solution $(A_x,A_y)$ in (1) of Theorem 7.3 
 can be chosen in such a way that they satisfy 
 the integrability condition $(5.13)$ at the initial value $(x_0,y)$, 
 \begin{equation}
 \frac{\partial A_y}{\partial x}(x_0,y)-\frac{\partial A_x}{\partial y}(x_0,y)
 +[A_x(x_0,y),A_y(x_0,y)]=0,
 \end{equation}
 for each $y$, 
 i.e., the initial functions 
 $f_0$, $f_1$ and $g_1$ may be chosen to satisfy that 
 \begin{equation}
 \frac{\partial A_x}{\partial y}(x_0,y)
 =g_1(y)+[f_0(y),f_1(y)]. 
 \end{equation}
\vskip0.6cm
Finally, we introduce a {\it loop group formulation} 
for biharmonic maps.
\par
We {\it first}, consider a ${\frak g}^{\mathbb C}$-valued 
$1$-forms 
\begin{equation}
\beta_{\nu}=\frac12(1-\nu)\,B_z\,dz+\frac12(1-\lambda^{-1})\,B_{\overline{z}}\,d\overline{z}
\end{equation}
for a parameter $\nu\in S^1$, 
which satisfy that 
\begin{equation}
d\beta_{\nu}
+[\beta_{\nu}\wedge\beta_{\nu}]=0
\quad(\forall\,\,\nu\in S^1), 
\end{equation}
where for 
the definition of $[\beta_{\nu}\wedge \beta_{\nu}]$, 
see (3.13). 
\par
{\it Next}, we consider 
${\frak g}^{\mathbb C}$-valued $1$-forms 
\begin{equation}
\alpha_{\nu}
=\frac12(1-\nu)\,A_z\,dz
+\frac12(1-\nu^{-1})\,A_{\overline{z}}\,d\overline{z}
\end{equation}
which satisfy that 
\begin{equation}
\left\{
\begin{aligned}
&\frac{\partial}{\partial z} 
\left(\delta\,\alpha_{\nu}\right)
-\left[\frac12(1-\nu)\,A_z,\delta\,\alpha_{\nu}\right]=B_z,\\
&\frac{\partial}{\partial \overline{z}} 
\left(\delta\,\alpha_{\nu}\right)
-\left[
\frac12(1-\nu)\,A_{\overline{z}},\delta\,
\alpha_{\nu}
\right]=
B_{\overline{z}},\\
&\quad d\,\alpha_{\nu}
+[\alpha_{\nu}\wedge\alpha_{\nu}
]=0,
\end{aligned}
\right.
\end{equation}
for each $\nu\in S^1$. 
Here, the co-differentiation 
$\delta\,\alpha_{\nu}$ of 
$\alpha_{\nu}$ is given by 
\begin{equation}
\delta\,\alpha_{\nu}
=-2\mu^{-2}\left(
\frac12(1-\nu)\,\frac{\partial}{\partial \overline{z}}
A_z
+\frac12(1-\nu^{-1})\,\frac{\partial}{\partial z}
A_{\overline{z}}
\right).
\end{equation}
Then, the mapping 
$\psi_{\nu}:\,\Omega\rightarrow G$
satisfying 
$\psi_{\nu}\,{}^{\ast}\theta=
\alpha_{\nu}$ is a biharmonic 
map $(\Omega,g)\rightarrow (G,h)$ where 
$g=\mu^{2}g_0$ for a positive $C^{\infty}$ 
function on $\Omega$.  
\vskip2cm\par       

\end{document}